\documentclass[10pt,leqno,twoside]{article}
\usepackage{amsmath,amssymb}
\author{Manuel Pedreira-P\'{e}rez
\and Luis-Eduardo Sol\'{a}-Conde\thanks{Supported by an F.P.U.
fellowship of Spanish Government.}}
\title{Projective generation and smoothness \\ of congruences of order
$1$\thanks{Mathematical Subject Classification: 14M15, 14J}}
\date{}

\newtheorem{teo}{Theorem}[section]
\newtheorem{defin}[teo]{Definition}
\newtheorem{propn}[teo]{Proposition}
\newtheorem{cor}[teo]{Corollary}
\newtheorem{lemma}[teo]{Lemma}
\newtheorem{example}[teo]{Example}
\newtheorem{rem}[teo]{Remark}
\def\rk{\mathop{\rm rank}}
\def\ord{\mathop{\rm ord}}
\def\cl{\mathop{\rm cl}}
\def\di{\mathop{\rm dim}}
\def\codi{\mathop{\rm codim}}
\def\mod{\mathop{\rm mod}}
\newcommand{\Cm}{\ensuremath{\textrm{%
\rm C\hspace{-0.98ex}\rule[0.3pt]{0.3pt}{6.5pt}\hspace{0.95ex}}}}
\newcommand{\Cn}{\ensuremath{\textrm{%
\rm C\hspace{-1.2ex}\rule[0.2pt]{0.3pt}{4.5pt}\hspace{0.95ex}}}}
\newcommand{\Pm}{\ensuremath{\textrm{\rm I\hspace{-1.8pt}P}}}
\def\qed{\hspace{\fill}$\rule{2mm}{2mm}$}
\newcommand\s{\sigma}
\newcommand\Si{\Sigma}
\newcommand\e{\emptyset}
\newcommand\op{\Omega}
\newcommand\lrw{\longrightarrow}
\newcommand\rw{\rightarrow}
\newcommand\mpt{\longmapsto}
\newcommand\T{\tau}
\newcommand\A{\alpha}
\newcommand\B{\beta}
\newcommand\G{\gamma}
\newcommand\la{\langle}
\newcommand\ra{\rangle}
\newcommand\I{{\cal I}}
\newcommand\Te{{\cal O}}
\font\euf=eufm10 at 12pt
\def\MA{\mbox{\euf m}}

\def\m@th{\mathsurround=\z@}

\begin{document}
\maketitle
\begin{abstract}
In this paper we give the projective generation of congruences of
order $1$ of $r$-dimensional projective spaces in $\Pm^N$
from their focal loci. In a natural way, this construction shows that the
corresponding surfaces in the grassmannian are the Veronese surface, and rational
ruled surfaces eventually with singularities. We characterize when these surfaces are smooth,
recovering and generalizing a Ziv Ran's result.
\end{abstract}
{\Large\bf Introduction}

\bigskip

Let $\Pm^N$ be the $N$--dimensional complex projective space 
and $G(r,N)$ be the Grassmannian of $\Pm^r$'s in $\Pm^N$. The study of surfaces
$\Si\subset G(r,N)$ --called {\it congruences}-- is a classical topic. The algebraic
equivalence class of $\Si$ consists of a pair of integers named {\it
order} and {\it class}: the order of the congruence is the number of elements of
$\Si$ meeting a general $(N-r-2)$--dimensional space. In 1866, Kummer
\cite{kummer} achieved the classification of congruences of lines in $\Pm^3$
of order $1$. A classification of congruences of order $1$ when $r\geq 1$ and $N\geq 3$ was
obtained independently by Ziv Ran in \cite{ran} using Bertini's theorem.  Kummer's basic idea
was studying the set of {\it Fundamental Points }of the congruence; that is, the points lying
on an infinite number of its lines. Recently F. Zak, A. Inshakov, S. Lvovski and A. Obolmkov
in
\cite{zak} have shown that the classification of congruences of order $1$ in $G(1,3)$ can be
rigorously done with elementary arguments. 

Since every fundamental point of a congruence is a focal point, one may ask whether the focal
methods provide sufficient information to classify congruences. The focal method was
introduced by Corrado Segre \cite{segre}, and it was used successfully in \cite{pedsol} to
obtain a birrational classification of congruences of planes; the modern reference here is
\cite{chiantini}. In this article we expose a classification of congruences of order $1$ for
all
$r\geq 1$ and $N\geq 3$ by using the focal method. The advantage of using it is that it provides
locally the construction of
$\Si$ from its focal locus $F(\Si)$. In particular, the method allows us to decide the
smoothness of the surface parametrizing it in the grassmannian recovering and generalizing
the Ziv Ran's result about smooth surfaces of order $1$ in $G(1,3)$. 

The paper is organized as follows:
In sections \ref{about} and \ref{focal} we summarize some properties
about surfaces in grassmannians and their focal loci. For a congruence of order $1$
we show that every focal point is a fundamental point. Another important ingredient here is the
relation between sections of a congruence in $\Pm^N$ with the projections of the corresponding
surface in the grassmannian, reducing some proofs to the case $r=1$. Excepting the family of
$\Pm^r$'s containing a fixed $\Pm^{r-1}$ in $\Pm^{r+2}$, the focal locus $F(\Si)$ has
dimension
$r$, obtaining three possibilities: either it is irreducible of degree $1$ (case III)
or $3$ (case I), or
$F(\Si)=\pi\cup X$, being $\pi$ a projective space and $X$ a rational
scroll (case II). It is also shown that if $F(\Si)$ contains a
projective space, $\Si$ is a rational ruled surface.
Section \ref{uno} deals with the case I: such a congruence is parametrized by a
Veronese surface; we also prove that there are only three fixed point free models, and the
remaining are constructed as cones over one of those.
In section \ref{dos}, we show the projective generation of a congruence $\Si$ in the
case II from $F(\Si)=\pi\cup X$, studying the normal models of the scroll $X$; we characterize  
the singular locus of $\Si$, and prove the existence of such surfaces for each invariant $e$.
In section \ref{tres} it is proved that a congruence in the case III is parametrized
by a rational cone, as well as its existence. We also give a characterization of
such congruences through a morphism between projective spaces.
Finally, in section \ref{lisas}, we resume all the results about smoothness of surfaces of order
$1$ in grassmannians.

\section{About Congruences of Order $1$}\label{about}
Let $\Pm^N$ be the $N$--dimensional projective complex space and $G(r,N)$
be the Grassmannian of $r$--dimensional subspaces in  
$\Pm^N$. We begin compiling some basic facts about $G(r,N)$ (see \cite{kleiman}
for more details):

It is well known that the cohomology group $H^i(G(r,N),{\mathbb{Z}})$ is $0$ if 
$i\notin[0,2(N-r)(r+1)]$ and the direct sum
$$
H^*(G(r,N),{\mathbb{Z}})=\bigoplus_i H^i(G(r,N),{\mathbb{Z}})
$$
is a graded ring with the cup--product $\cup$. We can assign the cohomology
class $[X]$ to each closed subset $X\subset G(r,N)$, being satisfied the next
properties:
\begin{itemize}
\item Two algebraically equivalent subvarieties have the same cohomology class.
\item If $X\subset G(r,N)$ is irreducible
$((r+1)(N-r)-p)$--dimensional, then $[X]\in H^{2p}(G(r,N),{\mathbb{Z}})$.
\item  For each two subvarieties $X,Y\subset G(r,N)$, $[X\cap Y]=[X]\cup[Y]$.
\end{itemize} 

Let $A_0\subsetneq A_1\subsetneq\ldots\subsetneq A_r$ be a strictly increasing sequence of
subspaces of
$\Pm^{N}$ , $a_i:=\di A_i$, 
$i=0,\ldots,r$, 
and $\op(A_0,\ldots,A_r)$ denote the next subvariety of the grassmannian:
$$
\op(A_0,\ldots,A_r)=\{\s\in G(r,N)\;/\di(\Pm^r(\s)\cap A_i)\geq
i,\,i=0,\ldots,r\}
$$
Such a subvariety is called a {\it Schubert cycle} in $G(r,N)$. Regarding
$G(r,N)$ as a projective variety through the Pl\"ucker's embedding, Schubert
cycles can be obtained cutting it with some linear subspaces; for instance, the
Schubert cycle $\op_r\subset G(1,3)$ parametrizing the lines in $\Pm^3$ meeting a
given one $r$ consists of the intersection $\op_r=G(1,3)\cap T_{G(1,3),[r]}$.

Given two Schubert cycles $\op(A_0,\ldots,A_r)$, $\op(B_0,\ldots,B_r)$, they
are algebraically equivalent when $\di A_i=\di B-i$ for every $i$; thus the
cohomology class of $\op(A_0,\ldots,A_r)$ depends only of the integers $a_i:=\di
A_i$, $i=0,\ldots,r$, and we denote it by $\op(a_0,\ldots,a_r)$. They play an
important role in the cohomology of $G(r,N)$, in fact:
\begin{teo}\label{basis}
$H^*(G(r,N),{\mathbb{Z}})$ is a free abelian group with basis
$$
\{\op(a_0,\ldots,a_r)\;/\,0\leq a_0<a_1<\ldots<a_r\leq N\}.
$$\qed
\end{teo}

Moreover, each $\op(A_0,\ldots,A_r)$ is irreducible of dimension
$\displaystyle\sum_{i=0}^r(a_i-i)$, so we have the next corollary:
\begin{cor}\label{bases} 
For every integer $p$, $H^{2p}(G(r,N),{\mathbb{Z}})$ is a free abelian group generated by 
the elements 
$\op(a_0,\ldots,a_r)$ with $p=[(r+1)(N-r)-\displaystyle\sum_{i=0}^r(a_i-i)]$.
Each
$H^{2p+1}(G(r,N),{\mathbb{Z}})$ is $0$. \qed
\end{cor}
\begin{example}\label{codos}
\begin{em}
By the above, we can calculate the basis of the cohomology groups of codimension
$2$ and dimension $2$:
$$
\begin{array}{l}
{H^{2(N-r)(r+1)-4}(G(r,N),{\mathbb{Z}})=\la\op_1:\equiv\op(0,1,\ldots
r-1,r+2),}\\
{\op_2:\equiv\op(0,1,\ldots r-2,r,r+1)\ra}\\
{H^{4}(G(r,N),{\mathbb{Z}})=\la\op^1:\equiv\op(N-r-2,N-r+1,N-r+2,\ldots
N),}\\
{\op^2:\equiv\op(N-r-1,N-r,N-r+2\ldots N)\ra}\\
\end{array}
$$
Denoting by $L^i\subset\Pm^n$ a fixed generic $i$--dimensional subspace, the
elements above correspond to the classes of the Schubert cycles:
$$
\begin{array}{l}
\op_1=[\{\s\in G(r,N)\;/L^{r-1}\subset\Pm(\s)\subset L^{r+2}\}]\\
\op_2=[\{\s\in G(r,N)\;/L^{r-2}\subset\Pm(\s)\subset L^{r+1}\}]\\
\op^1=[\{\s\in G(r,N)\;/\Pm(\s)\cap L^{N-r-2}\neq\e\}]\\
\op^2=[\{\s\in G(r,N)\;/\di(\Pm(\s)\cap L^{N-r})\geq 1\}]\\
\end{array}
$$
and we can compute the products $\op_i\cup\op^j$ as the classes of the
intersections of such Schubert cycles; thus $\op_i\cup\op^j=0$ if $i\neq j$ and
$\op_i\cup\op^j=1$ if $i=j$. Hence the cohomology class of a $2$--dimensional
irreducible subvariety $\Si\subset G(r,N)$ is determined by the number of points
of intersection of $\Si$ with a generic Schubert cycle in the class of $\op^1$
(the number of elements of $\Si$ cutting a generic $(N-r-2)$--dimensional space),
and with another one in the class of
$\op_2$ (the number of elements of $\Si$ meeting a generic $(N-r)$--dimensional
space in lines); these numbers are called, respectively, {\it order}
($\ord(\Si)$) and {\it class} ($\cl(\Si)$) of $\Si$. The pair
($\ord(\Si,\cl(\Si)$) is called {\it bidegree} of $\Si$. Moreover, regarding
$\Si$ as a projective variety through the Pl\"ucker embedding of the
grassmannian, the degree of $\Si$ is exactly $\ord(\Si)+\cl(\Si)$.
\end{em}
\end{example}

We wish to investigate $2$--dimensional families of $r$--dimensional subspaces
in $\Pm^N$, called classically {\it congruences}. From now on,
$\Si$ will denote a
$2$--dimensional reduced and irreducible closed subvariety of $G(r,N)$ 
and $V(\Si)$ its projective realization, 
\begin{equation}
V(\Si):=\bigcup_{\s\in\Si}\Pm^r(\s)
\label{realiza1}
\end{equation}

The objective of this work is the classification of congruences of order $1$;
by the next lemma we will be reduced to studying the case $N=r+2$.
\begin{lemma}\label{r+2}
Let $\Si\subset G(r,N) $ be a congruence. If $\ord(\Si)=1$, then there exist
a $(r+2)$--dimensional projective subspace {\rm $A\subset \Pm^N$} such that 
{\rm $\Pm^r(\s)\subset A$} for all $\s\in \Si$.
\end{lemma}

{\bf Proof:} Since $\Si$ is irreducible, $V(\Si)$ is irreducible too. Let
$B\subset\Pm^N$ be a generic $(N-r-2)$--dimensional subspace. As $\ord(\Si)=1$ we have
$B\cap V(\Si) =B\cap\Pm^r(\s)\neq\e$
for only one $\s\in\Si$. It follows that $\di V(\Si)=r+2$ and $\deg V(\Si)=1$. We
finish taking $A:=V(\Si)\simeq \Pm^{r+2}$. \qed

According to this lemma, we will only consider $\Si\subset G(r,r+2)$, using the
next notation for the incidence variety:
\begin{center}
\setlength{\unitlength}{5mm}
\begin{picture}(18,3)
\put(8,2.7){\makebox(0,0){$\supset{\cal
I}_{\Sigma}:=\,\{(P,\sigma)\in {\bf
P}^{r+2}\times\Sigma;\,P\in {\bf P}^r(\sigma)\}$}}
\put(0,2.7){\makebox(0,0){${\bf
P}^{r+2}\times\Sigma$}} 
\put(14.8,2.7){\vector(1,0){2}}
\put(2.8,2.1){\vector(0,-1){1.5}}
\put(3.2,1.5){\makebox(0,0){$q$}}
\put(16,3.1){\makebox(0,0){$p$}}
\put(18,2.7){\makebox(0,0){${\bf P}^{r+2}$}}
\put(2.8,0){\makebox(0,0){$\Sigma$}}
\end{picture}
\end{center}
\begin{defin}\label{funda}
{\em Let $\Si\subset G(r,r+2)$ be a congruence. A point $P\in\Pm^{r+2}$ is called
{\it fundamental} for $\Si$ if
$\di(p^{-1}(P))\geq 1$. The congruence $\Si$ is called {\it degenerate} iff
$\di V(\Si)<r+2$, that is iff every point in $V(\Si)$ is fundamental. For
instance, a congruence of order $1$ is nondegenerate.}
\end{defin}

Given a congruence $\Si$ in $\Pm^{r+2}$ and a projective subspace
$L\subset\Pm^{r+2}$, we can obtain a congruence in $L$ cutting it with each
element of $\Si$. The next results provide some properties of such sections.
\begin{lemma}\label{cortar}
Let $\Si\subset G(r,r+2)$ be a congruence of 
order $\ord(\Si)$ and $L\subset\Pm^{r+2}$ be a $(l+2)$--dimensional
generic subspace, $l\geq 1$; then $\Si$ does not meet $\op_L:=\{\s\in
G(r,r+2)\;/\di L\cap\Pm^r(\s)\geq l+1\}$, which is a codimension $l+2$ Schubert 
cycle; thus the projection
$$\rho_L:\Si\subset G(r,r+2)\lrw
G(l,L),\quad\pi_L(\s):=[\Pm^r(\s)\cap L]
$$
 is regular and has as image a surface $\Si_L$ such that
$\ord(\Si_L)=\ord(\Si)$, $\cl(\Si_L)=\cl(\Si)$.
\end{lemma}
{\bf Proof:} Consider the incidence variety
$$
\I:=\{(L,\sigma)\,/\di(L\cap\Pm^r(\s))\geq
l+1\}\subset G(l+2,r+2)\times\Si
$$
with projections $p_1$ and $q_1$. Given $\s\in\Si$, $q_1^{-1}(\s)\equiv\{L\in
G(l+2,r+2)\;/\di(\Pm^r(\s)\cap L)\geq l+1\}$ is a
codimension $l+2$ Schubert cycle (because its cohomology class is
$\op(r-1-l,\ldots,r-1,r,r+2)$), so
$\di \I=\di\Si+\di q_1^{-1}(\s)=2+\di G(l+2,\Pm^{r+2})-l-2=\di
G(l+2,\Pm^{r+2})-1$. Therefore 
$p_1(\I)\subset G(l+2,r+2)$ is a proper closed subset, and if $L\in
G(l+2,r+2)\setminus p_1(\I)$, then $\di(L\cap\Pm^r(\s))=l$ for all $\s\in\Si$. It
follows that the map $\rho_L$ defined above is regular. Furthermore it is
birrational: a generic point $P\in\Pm^{r+2}$ is contained in a finite number of
elements of the congruence,
$\Pm^r(\s_1),\ldots,\Pm^r(\s_n)$; a generic $l$--dimensional subspace 
through $P$, $W\subset
\Pm^r(\s_1)$ will not be contained in another $\Pm^r(\s_i)$; so taking
$L\cap\Pm^r(\s_1)=:W$, we will have 
$L\cap\Pm^r(\s_1)=W\neq
L\cap\Pm^r(\s_2),\ldots,L\cap\Pm^r(\s_n)$; thus
$\rho_L^{-1}(W)=\{\s_1\}$, so $\deg\rho_L=1$.

By construction, the equality $\ord(\Si_L)=\ord(\Si)$ holds; by
the birrationality of $\rho_L$ $\deg(\Si_L)=\deg(\Si)$ and so 
$\cl(\Si_L)=\cl(\Si)$ holds too. \qed
\begin{rem}\label{proyeccion}
{\em Furthermore, regarding $G(r,r+2)$ and
$G(l,l+2)$ as projective varieties through the corresponding Pl\"ucker
embeddings, the map $\rho_L$ is easily shown to be the projection from the
space generated by the Schubert cycle 
$\{\s\in G(r,r+2)\;/\di(\Pm^r(\s)\cap L)\geq 2\}$.}
\end{rem}
\begin{lemma}\label{dimfunda}
Let $\Si\subset G(r,r+2)$ be a nondegenerate congruence, and
$F(\Si)\subset\Pm^{r+2}$ its fundamental locus. Then $\di
F(\Si)\leq r$. 
\end{lemma}
{\bf Proof:} Suppose the lemma were false. We could find an
$(r+1)$--dimensional irreducible component $C\subset F(\Si)$.
Considering the restriction $p^{-1}(C)\subset\I_{\Si}\lrw
C\subset\Pm^{r+2}$, we would have $\di p^{-1}(C)=r+2$, and so
$p^{-1}(C)=\I_{\Si}$. Consequently $\Si$ would be degenerate.
\qed
\begin{lemma}\label{cortfunda}
Under the hypotheses of (\ref{cortar}) and 
with the above notation, it is verified the equality
$F(\Si)\cap L=F(\Si_L)$.
\end{lemma}
{\bf Proof:} By construction $F(\Si)\cap L\supset F(\Si_L)$. If $P\in
F(\Si)\cap L\setminus F(\Si_L)$, $P$ would be contained in an
irreducible family of elements of $\Si$ with the same trace on $L$. An
$l$--dimensional subspace contained in infinite elements of $\Si$ is called {\it
fundamental}. Suppose
$\Si$ is nondegenerate (the degenerate case is trivial). It
is sufficient to show that, if
$D$ is an irreducible component of $F(\Si)$ and if $L$ is generic, $D\cap L$
does not contain fundamental subspaces. On the contrary, suppose $D\cap
L$ contains fundamental subspaces for all $L$; by 
(\ref{dimfunda}) $\di D\leq r$, so $D$ would
be an $r$--dimensional projective space whose $l$--dimensional subspaces are all
fundamental. Considering the incidence variety: 
$$
{\cal J}:=\,\{(l,\sigma)\,/l\subset\Pm^r(\s)\cap D\}\subset G(l,D)\times\Si
$$
with projections $p_2$ and $q_2$,
we have $\di{\cal J}=2+\di q_2^{-1}(\s)\leq 2+(l+1)(r-1-l)$,
because $q_2^{-1}(\s)\cong G(1,\Pm^r(\s)\cap D)$. Thus $\di{\cal
J}\leq 2+(l+1)(r-1-l)\leq(l+1)(r-l)=\di G(l,D)$, so either $p_2$ is not
surjective, or it has generically finite fibers. In both cases, the generic $l$--dimensional
subspace is not fundamental, a contradiction. \qed
\begin{defin}\label{pfijosdef}
{\em  A point $P\in\Pm^r(\s)$ is called {\it fixed
point} of the congruence $\Si$ when $P\in\Pm^r(\s)$ for
all $\s\in\Si$.}
\end{defin}
\begin{rem}\label{greater}
{\em The locus of fixed points of a congruence $\Si$
is a linear subspace that we will denote by $T(\Si)\subset\Pm^{r+2}$, and
$\Si$ lies on the Schubert's cycle:
$$\op_{T(\Si)}:=\{\s\in G(r,r+2)/\Pm^r(\s)\supset T(\Si)\}$$ Being
$L\subset\Pm^{r+2}$ a generic complementary subspace of $T(\Si)$, $k+2:=\di
L=r+2-\di T(\Si)-1$, the maps
\begin{equation}
\begin{array}{c}
\rho_L:\s\in\op_{T(\Si)}\mpt[\Pm^r(\s)\cap L]\in G(k,k+2)\\
\rho^{T(\Si)}:\T\in G(k,k+2)\mpt[\Pm^k(\T)+T(\Si)]\in\op_{T(\Si)}\\
\end{array}
\label{prefijos1}
\end{equation} 
provide an isomorphism $G(k,k+2)\cong\op_{T(\Si)}\subset G(r,r+2)$. Clearly,
the congruence $\rho_L(\Si)\cong\Si$ has no fixed points; roughly speaking,
$\rho_L(\Si)$ is the same surface as $\Si$ living in a grassmannian of lesser
dimension. Therefore {\it it suffices to study congruences without fixed
points}}.
\end{rem}
\section{Focal Locus of a Congruence in $G(r,r+2)$}\label{focal}
\begin{defin}\label{focaldefin}
{\em Let $\Si\subset G(r,r+2)$ be a congruence. Being $\s\in\Si$ smooth, a point
$P\in\Pm^r(\s)$ is called {\it focal} if $(dp)_{(P,\s)}$ is not injective. If
$\ker((dp)_{(P,\s)})\supset\la v\ra$ for some $v\in q^*T_{\Si,\s}$,
$v\neq 0$, we say that $P$ is focal for the direction $\la v\ra$.}
\end{defin}

\begin{propn}\label{equiv}
Let $\Si\subset G(r,r+2)$ be a congruence. If
$P\in\Pm^{r+2}$ is a fundamental point of $\Si$, then $P$ is focal for
every $\s\in\Si$ smooth with $P\in\Pm^r(\s)$. Furthermore, if
$P\in\Pm^r(\s)$ is focal and if $\ord(\Si)=1$, then $P$ is fundamental.
\end{propn}
{\bf Proof:} Being $P$ fundamental, there is an irreducible curve
$C\in\Si$ such that $P\in\s$ for all $\s\in C$. Taking a smooth point
$\s\in C$ and regarding the curve $\{P\}\times C\subset\I_{\Si}$, we
have $p(\{P\}\times C)=P$. Assuming that $\s$ is smooth in $\Si$ we can
write $(dp)_{(P,\s)}(T_{\{P\}\times C,(P,s)})=0$. Therefore
$P\in\Pm^r(\s)$ is focal.

For the second part we use that $\ord(\Si)=\deg(p)$. If
$\ord(\Si)=1$, the map $p$ is dominant and $p(\I_{\Si})=\Pm^{r+2}$ is
normal. Hence $\deg(p)=\sum_{p(x)=P}m_x(p)$ for every nonfundamental
point $P$, being
$m_x(p):=\di_{\Cn}(\Te_{\I_{\Si},x}/p^*\MA_P)$ for $\MA_P$ the maximal ideal
of the local ring $\Te_{\Pm^{r+2},P}$. Thus, if $\ord(\Si)=1$,
$m_x(p)>1$ forces $P$ to be fundamental. The proof is completed with
the next lemma.\qed
\begin{lemma}\label{multi}
Let $\Si\subset G(r,r+2)$ be a congruence and denote
$x:=(P,\s)\in\I_{\Si}$. $x$ is focal if and only if $m_x(p)\geq  2$.
\end{lemma}   
{\bf Proof:} Regarding the diagram
\begin{center}
\setlength{\unitlength}{5mm}
\begin{picture}(16,3)
\put(1,0){\makebox(0,0){$\MA_x$}}
\put(1,3){\makebox(0,0){$p^*\MA_P$}}
\put(8,0){\makebox(0,0){$\Te_{\I_{\Si},x}$}}
\put(8,3){\makebox(0,0){$\Te_{\I_{\Si},x}$}}
\put(15,0){\makebox(0,0){$\Cm$}}
\put(15,3){\makebox(0,0){$\Te_{\I_{\Si},x}/p^*\MA_P$}}
\put(1,2.2){\vector(0,-1){1.5}}
\put(0.85,2){\line(1,0){0.3}}
\put(7.9,2.2){\line(0,-1){1.5}}
\put(8.1,2.2){\line(0,-1){1.5}}
\put(15,2.2){\vector(0,-1){1.5}}
\put(14.85,1.1){\line(1,0){0.3}}
\put(3,0){\vector(1,0){3}}
\put(3,3){\vector(1,0){3}}
\put(10,0){\vector(1,0){2.8}}
\put(10,3){\vector(1,0){2.8}}
\end{picture}
\end{center}
it is easily seen that 
$\di_{\Cn}(\Te_{\I_{\Si},x}/p^*\MA_P)=\di_{\Cn}(\MA_x/p^*\MA_P)+1$. Consider the
next diagram:
\begin{center}
\setlength{\unitlength}{5mm}
\begin{picture}(16.1,3)
\put(1,0){\makebox(0,0){$\MA_x^2$}}
\put(1,3){\makebox(0,0){$\MA_P^2$}}
\put(6,0){\makebox(0,0){$\MA_x$}}
\put(6,3){\makebox(0,0){$\MA_P$}}
\put(13.4,0){\makebox(0,0){$\MA_x/\MA_x^2=(T_{\I_{\Si},x})^*$}}
\put(13.7,3){\makebox(0,0){$\MA_P/\MA_P^2=(T_{\Pm^{r+2},P})^*$}}
\put(1,2.2){\vector(0,-1){1.5}}
\put(0.85,2){\line(1,0){0.3}}
\put(6,2.2){\vector(0,-1){1.5}}
\put(5.85,2){\line(1,0){0.3}}
\put(15,2.2){\vector(0,-1){1.5}}
\put(2,0){\vector(1,0){3}}
\put(2,3){\vector(1,0){3}}
\put(7,0){\vector(1,0){2.8}}
\put(7,3){\vector(1,0){2.8}}
\put(1.5,1.5){\makebox(0,0){$p^*$}}
\put(6.5,1.5){\makebox(0,0){$p^*$}}
\put(16.1,1.5){\makebox(0,0){$(dp)_x^*$}}
\end{picture}
\end{center}
By {\it Nakayama's lemma} $p^*$ is surjective iff $(dp)_x^*$
is surjective, that is iff $(dp)_x$ is injective. Hence, $x$
is focal iff $\di_{\Cn}(\MA_x/p^*\MA_P)>0\iff m_x(p)>1$. \qed

We continue with some basic facts about the focal locus of a congruence in
$G(r,r+2)$; we refer the reader to \cite{pedsol} for more details in the
case $r=2$, but the proofs also work for all $r$. Let
$\Si\subset G(r,r+2)$ be a congruence and
$\s\in\Si$ an smooth point. For each $P\in\Pm^r(\s)$ there is an
isomorphism $T_{\I_{\Si},(P,\s)}\cong T_{\Si,\s}\oplus T_{\Pm^r(\s),P}$.

The {\it focal locus at} $\Pm^r(\s)$ can be
calculated through the {\it Characteristic Map of
Kodaira--Segre--Spencer}:
\begin{equation}
\chi:T_{\Si,\s}\quad\lrw\quad
H^0({\cal N}_{\Pm^r(\s),\Pm^{r+2}})\cong H^0({\cal
O}_{\Pm^r(\s)}(1))\oplus H^0({\cal O}_{\Pm^r(\s)}(1))
\label{caesese}
\end{equation}
being $\chi(v)=0$ the equations of the focal locus of $\Si$ at
$\Pm^r(\s)$ in the direction $\la v\ra\subset T_{\Si,\s}$. An easy
computation shows that:
\begin{itemize}
\item The focal locus at $\Pm^r(\s)$ in a direction $\la
v\ra$ is a projective subspace of dimension $r-2$ or $r-1$
and we denote it by $Z(\chi(v))$. Taking a base $\{v_1,v_2\}$ in 
$T_{\Si,\s}$, whose images through $\chi$ are $(f_{11},f_{12})$ and 
$(f_{21},f_{22})$ respectively, $Z(\chi(\lambda v_1+\mu v_2))$ will be the
set of points $P\in\Pm^r(\s)$ satisfying the equations:
$$
\left.\begin{array}{c}
{\lambda f_{11}(P)+\mu f_{21}(P)=0}\\
{\lambda f_{12}(P)+\mu f_{22}(P)=0}\\
\end{array}\right\}
$$
\item The focal locus at $\Pm^r(\s)$ is a quadric  $Q(\s)=\bigcup_{v\in
T_{\Si,\s}}Z(\chi(\s))$, and it is given by the equation 
\begin{equation}
\det\left(\begin{array}{cc}
{f_{11}(P)}&{f_{21}(P)}\\
{f_{12}(P)}&{f_{22}(P)}\\
\end{array}\right)=0
\label{cuadfocal1}
\end{equation}
It is reducible
if and only if there exist a direction
$\la v\ra\subset T_{\Si,\s}$ such that $\di Z(\chi(v))=r-1$. Such
directions are called {\it developable}, and the number of
them is fixed in an open subset of $\Si$; this number can be
$0$, $1$, $2$ (distinct or coincident) or $\infty$ (if all
direction is developable).
\item If $Q(\s)$ is irreducible, it contains a
$1$--dimensional family of $(r-2)$--dimensional projective
subspaces. Therefore $\rk Q(\s)\leq 4$ (see for instance
\cite{enriques}, book V, page 100).
$Q(\s)$ will be a cone with vertex a subspace
$C(\s)\subset\Pm^r(\s)$ over an smooth conic or an smooth
quadric surface. In both cases the focal loci in each
direction are subspaces of maximal dimension contained in
$Q(\s)$, and so the intersection of two of them is $C(\s)$. Hence 
$C(\s)$ is the set of points focal for all direction.
\item Let $\Delta\subset G(r,r+2)$  be a reduced and irreducible
curve. We say that $\Delta$ is a {\it developable system} if
the focal locus in each $\Pm^r(\s)$,  $\s\in\Delta$ has dimension
$r-1$. Following 
\cite{pedsol}, $\Delta$ parametrizes a family of the next
kind: Cone with vertex a subspace $\Pm^k$ over the family of
$(r-k-1)$--osculating spaces to an irreducible curve
${\cal C}\subset\Pm^{r-k+1}$, $k=-1,0,\ldots,r-1$. For instance, a
$1$--dimensional family of $\Pm^r$'s lying in an $(r+1)$--dimensional space
is always developable.
\item A nondegenerate congruence $\Si\subset G(r,r+2)$
has reducible generic focal quadric iff by the
generic $\s\in\Si$:
\begin{itemize}
\item is passing $1$ developable system, or
\item are passing $2$ distinct developable systems, or
\item are passing $2$ coincident developable systems, or
\item are passing infinite developable systems. 
\end{itemize}
\end{itemize}

Being $\Si\subset G(r,r+2)$, denote
$F_1(\Si)\subset\Pm^{r+2}$ its focal variety, that is, the projective
realization of the family of focal quadrics, which is defined in an open set
$U\subset\Si$, ${\cal F}_1(U)\lrw U$.  In general, the fundamental locus
$F(\Si)$ is a subvariety of $F_1(\Si)$ (see \cite{ciliberto}). In the case 
$\ord(\Si)=1$, (\ref{equiv}) clearly forces the equality 
$F_1(\Si)=F(\Si)$. Moreover, by (\ref{dimfunda}) $\di F(\Si)\leq r$, 
$\di F(\Si)<\di F_1(U)=r+1$. The generic focal quadric can be irreducible or
not, so either $F(\Si)$ is irreducible, or it is the projective realization of
two $2$--dimensional families of $\Pm^{r-1}$'s, having at least $2$ irreducible
components.
\begin{example}\label{c1c2}{\em Let $C_1$, $C_2$ be two irreducible curves in $\Pm^3$.
We write $\Si(C_1,C_2)$ for the set of lines joining points of both curves, and $\Si(C_1)$ the
family of secant lines to $C_1$. Clearly the focal locus of these congruences
contains the base curves (furthermore, they will be exactly the focal locus in some
cases). If
$C_1$ and
$C_2$ are not coplanar, then
$\Si(C_1,C_2)$ is nondegenerate, and if $C_1$ is not a plane curve, then
$\Si(C_1)$ is nondegenerate. We will show later when these
congruences have order $1$.} 
\end{example}
In the next two results we examine the focal locus of a congruence of order $1$,
concretely its irreducible components and degree.
\begin{propn}\label{compo}
If the focal locus of a congruence $\Si$ of order $1$ is reducible,
$F(\Si)=C_1\cup C_2$, then both components have dimension $r$ and one of
them is linear.
\end{propn}
{\bf Proof:} We have shown that the number of irreducible components of $F(\Si)$ is
at most two, with dimensions lesser than or equal to $r$ by (\ref{dimfunda}). Let
$C_1$, $C_2$ be the irreducible components of $F(\Si)$. It is sufficient to show that
if
$F(\Si)$ is reducible, its section with the generic
$3$--dimensional space consists of a line and an irreducible curve. Thus using 
(\ref{cortar}) and
(\ref{cortfunda}), we can suppose
$r=1$. If
$\di C_i=0$ for some $i$, then
$\Si$ would be the family of lines passing by $C_i$, whose focal locus is the point
$C_i$, so $C_1$ and $C_2$ are irreducible curves. Every line of $\Si$ contains a
point in each curve, so $\Si=\Si(C_1,C_2)$, and $\Si\cap\Si(C_1)$ and $\Si\cap\Si(C_2)$
are closed subsets in $\Si$ (if $\Si=\Si(C_1)$, then every line of $\Si$ would
contain at least three focal points, two of $C_1$ and one of $C_2$, so the focal
quadric at each line must be the whole line; equivalently, the congruence would be
degenerate, so $\ord(\Si)=0$, false.). Being
$K\subset\Pm^3$ a generic plane, it contains
$\cl(\Si)$ generators of $\Si$, none belonging to $\Si(C_1)$ nor $\Si(C_2)$. Hence  
every generator in $K$ meets $C_1$ and $C_2$ exactly one time. If $\deg
C_1,\deg C_2\geq 2$, we could take $l=P_1P_2$, $l'=Q_1Q_2$, $P_i,Q_i\in C_i$,
$P_i\neq Q_i$. The point $P=l\cap l'$ could not lie on $C_1$, $C_2$,
but as $\ord(\Si)=1$, $P$ should be fundamental, a contradiction. \qed
\begin{propn}\label{irred}
If the focal locus $F(\Si)$ of a congruence of order $1$ is irreducible, one of these
possibilities is true:
\begin{enumerate}
\item $F(\Si)$ is an $(r-1)$--dimensional projective space and
$\Si$ is the family of $\Pm^r$'s containing $F(\Si)$.
\item $\di F(\Si)=r$, and then:
\begin{enumerate}
\item either $F(\Si)$ is a projective space,
\item or $F(\Si)$ is an irreducible variety of degree $3$.
\end{enumerate}
\end{enumerate}
\end{propn}
{\bf Proof:} By hypothesis $F(\Si)$ is irreducible, containing the focal quadrics. If
$\di F(\Si)\leq r-1$, then $F(\Si)=Q(\s)\subset\Pm^r(\s)$ for all $\s\in\Si$. Therefore
$F(\Si)$ is an  
$(r-1)$--dimensional projective space.

Suppose $\di F(\Si)=r$. Using (\ref{cortar}) and
(\ref{cortfunda}) we are reduced to the  case $r=1$. If $\deg F(\Si)=n>1$, then
$\Si=\Si(F(\Si))$: let
$\s\in\Si$ be generic and $P\in F(\Si)\cap\Pm^1(\s)\neq\e$; since $n>1$, a generic plane 
$\pi\supset\Pm^1(\s)$ contains another point $Q\in F(\Si)$, $Q\neq P$. A $1$--dimensional
family of lines of $\Si$ is passing by $Q$, forming a cone of
degree $\geq 1$; one of its generators $\Pm^1(\T)$ lies on $\pi$, cutting
$\Pm^1(\s)$ in a fundamental point $R\in F(\Si)$, and $\Pm^1(\s)=PR$.

Since $\Si=\Si(F(\Si))$, necessarily $\deg F(\Si)\geq 3$ (otherwise, $F(\Si)$ would
be planar, and $\Si$ would be degenerate). Suppose, contrary to our claim,
that $\deg F(\Si)\geq 4$. The generic plane $\pi\subset\Pm^3$ cuts $F(\Si)$ in
$n\geq 4$ points such that any $3$ of them are not collinear (by the Trisecant
Lemma). Taking four of them, the diagonal points of the square they form must be
fundamental; that is, they belong to $F(\Si)$; hence there are $3$ collinear points of
$F(\Si)\cap\pi$, a contradiction. \qed

We study now the case where $F(\Si)$ (irreducible or not) contains a projective
space.
\begin{lemma}\label{cone}
Let $\pi$ be an $r$--dimensional projective space and consider the Schubert's
cycle:
\begin{equation}
\op_{\pi}:=\{\s\in G(r,r+2)\,/\,\di\Pm^r(\s)\cap\pi\geq r-1\}\subset G(r,r+2)
\label{omegapi1}
\end{equation}
$\op_{\pi}$ is a cone with vertex $[\pi]$ over a Segre Variety of the form
$\Pm^r\times\Pm^1\subset\Pm^{2r+1}$
\end{lemma}  
{\bf Proof:} Regard $H_{\pi}=\{[H]\in{\Pm^{r+2}}^*\,/\,H\supset\pi\}\cong\Pm^1$
and the projection
\begin{equation}
\phi:\op_{\pi}\setminus\{[\pi]\}\lrw\pi^*\times
H_{\pi},\quad\phi(\s):=(\Pm^r(\s)\cap\pi,\Pm^r(\s)+\pi)
\label{fifi}
\end{equation}
whose fibers are $\phi^{-1}(h,H)=\op_{h,H}=\{\s\in
G(r,r+2)\,/\,h\subset\Pm^r(\s)\subset H\}\cong\Pm^1$, and $[\pi]\in\op_{h,H}$.
In this way $\phi$ defines $\op_{\pi}$ as a cone. \qed
\begin{rem}\label{conedos}
{\em Segre variety $\pi^*\times H_{\pi}$ can be geometrically thought in the next
way: a line $l\subset\Pm^{r+2}$ disjoint of $\pi$ parametrizes the elements
of $H_{\pi}$, so the subvariety of the grassmannian
$\op_{\pi}\cap\{\s/\,\Pm^r(\s)\cap l\neq\e\}$ is a Segre variety of the form 
$\pi^*\times H_{\pi}$. Moreover, it is the intersection of two Schubert
cycles, one of them not containing $[\pi]$, so it is a hyperplane section of
$\op_{\pi}$.}
\end{rem}
\begin{propn}\label{ratrul}
Let $\Si\subset G(r,r+2)$  be a congruence and  
$\pi\subset\Pm^{r+2}$ be an $r$--dimensional projective subspace such that
$\Si\subset\op_{\pi}$. $\Si$ has order $1$ if and only if
$\phi^{-1}(\pi^*\times\{H\})\cap\Si$ is a line for all $H\in H_{\pi}$. Thus the map
\begin{equation}
\A:\Si\setminus\{[\pi]\}\lrw H_{\pi}\cong\Pm^1
\label{alfalfa}
\end{equation}
given by $\A(\s)=\Pm^r(\s)+\pi$
defines $\Si$ as a rational ruled surface.
\end{propn}
{\bf Proof:} Taking $H\in H_{\pi}$, $\A^{-1}(H)=\{\s\in\Si\,/\,\Pm^r(\s)\subset
H\}$ is $1$--dimensional and each one of its irreducible components is a
developable family of $\Pm^r$'s (because they are contained in $H$), whose
focal locus has dimension
$\leq r$. If
$\ord(\Si)=1$, then $\di F(\Si)=r$, and so the focal locus of those families must be
$(r-1)$--dimensional. Equivalently, each one is  a family of $\Pm^r(\s)$'s containing
an $(r-1)$--dimensional space in $H$; since
$\ord(\Si)=1$, each $H\supset\pi$ can only contain one of such families
(otherwise, every point in $H$ would be focal). So $\A^{-1}(H)\cong\Pm^1$.
Conversely, a congruence constructed in such way has clearly order $1$. \qed
\begin{propn}\label{fixed}
Let $\Si\subset\op_{\pi}\subset G(r,r+2)$ be a congruence of
order $1$. If $\Si$ has no fixed points, then $\cl(\Si)\geq r$.
\end{propn}
{\bf Proof:} Consider
\begin{equation}
\B:\Si\setminus\{[\pi]\}\stackrel{\phi}{\lrw}\pi^*\times
H_{\pi}\stackrel{p}{\lrw}\pi^*
\label{betata}
\end{equation}
given by $\B(\s)=\Pm^r(\s)\cap\pi$. In each hyperplane $H\supset\pi$, $\Si$ consists
of a pencil of $\Pm^r$'s contained in $H$ (with base a certain
$(r-1)$--dimensional projective space). If
$P\in\Pm^{r+2}$ is a fixed point of $\Si$, necessarily $P\in\pi$, and so
$\B(\Si)\subset\pi^*(P)=\{h\in\pi^*\,/P\in h\}\subset\pi^*$ will be degenerate. 
Thus $\Si$ has fixed points iff $\B(\Si)\subset\pi^*$ is degenerate.
Now, we have two cases:
\begin{enumerate}
\item If $\phi(\Si)$ is a curve: 
$\cl(\Si)+1=\deg(\Si)=\deg(\phi(\Si))=\#(\phi(\Si)\cap(h\times
H_{\pi}+\pi^*\times H))$ where $h\in\pi^*$ and $H\in H_{\pi}$. Hence
$\cl(\Si)+1=\#(\phi(\Si)\cap(h\times H_{\pi}))+\#(\phi(\Si)\cap(\pi^*\times
H))\geq\deg(\B(\Si))+1\geq\di(\B(\Si))+1=r+1$.
\item If $\phi(\Si)$ is a surface: $\cl(\Si)+1=\deg(\Si)\geq\deg(\phi(\Si))$. The
map $p:\phi(\Si)\lrw\B(\Si)$ consists of projecting from a certain $\pi^*\times
H$, that contains exactly one generator of $\phi(\Si)$. Hence
$\cl(\Si)+1\geq\deg(\B(\Si))+2\geq\di\la\B(\Si)\ra-1+2=r+1$. \qed 
\end{enumerate}
 \section{Case I: $F(\Si)$ is $r$--dimensional irreducible of degree
$3$}\label{uno}  First examine the case $r=1$: by (\ref{irred})
$F(\Si)$ is a nondegenerate curve of degree $3$ in $\Pm^3$, that is a rational normal
cubic, and
$\Si=\Si(F(\Si))$. Conversely, the congruence of secant lines to a rational normal
cubic has clearly order $1$. Since a generic plane
$\pi\subset\Pm^3$ cuts $F(\Si)$ in three noncollinear points, $\cl(\Si)=3$, thus 
$\deg(\Si)=4$.
Moreover, the secant lines through a fixed point $P\in F(\Si)$ form a quadric cone, so 
$\Si$ is not a
ruled surface (the lines in $G(1,3)$ parametrize linear pencils of lines through a point).
Such a quadric cone is parametrized by an irreducible conic $\Si(P)\subset\op_P$
where $\op_P$ is the family of lines through $P$. Given two points
$P,Q\in F(\Si)$, $\op_P\cap\op_Q=[PQ]$ so
$\la\Si(P),\Si(Q)\ra=\op_P+\op_Q=T_{G(1,3),[PQ]}$ that is $4$--dimensional
and $T_{G(1,3),[PQ]}\cap\Si=\Si(P)\cup\Si(Q)$. Hence 
$\Si$ is nondegenerate in $\Pm^5=\la G(1,3)\ra$. By Del Pezzo's
Theorem (see for instance \cite{griffiths} page 525),
$\Si$ must be the Veronese surface.


By (\ref{cortar}), if $r>1$, a generic $3$--dimensional subspace
$L\subset\Pm^{r+2}$ provide us a projection $\rho_L:\Si\lrw\Si_L\subset\Pm^5$, where
$\Si_L$ is a Veronese surface by the above. Since $\Si_L$ is normal and $\rho_L$
is regular and birrational, $\Si$ will be a Veronese surface too.
 
$F(\Si)$ is an $r$--dimensional irreducible variety of degree $3$ in $\Pm^{r+2}$. A
description of such varieties can be found in \cite{eisen}: Such a variety is, in general, a
cone with vertex a subspace $C$ over a variety of the same degree and codimension in a
projective subspace complementary to $C$. If it is smooth, then it is a rational normal
scroll. Thus it contains a
$1$--dimensional family of disjoint $\Pm^{r-1}$'s, which is only possible if
$r\leq 3$.

Therefore $F(\Si)$ is one of the next varieties:
\begin{itemize}
\item the rational normal cubic in
$\Pm^3$ ($r=1$),
\item the rational normal cubic surface in $\Pm^4$, projective realization of the Blowing-Up
of the plane in a point ($r=2$),
\item the rational normal cubic $3$--fold in $\Pm^5$, projection
of the Blowing-Up of $\Pm^3$ in a line from one of its points ($r=3$), 
\item a cone with vertex a subspace $C$ over one of the varieties described above.
\end{itemize}
Such varieties
are intersection of
$3$ quadrics, and they contain exactly a $2$--dimensional family of $(r-1)$--dimensional
quadrics. These quadrics are smooth iff $F(\Si)$ is smooth, othercase they are cones with vertex
$C$. $\,\Si$ must therefore be the congruence of $\Pm^r$'s containing such quadrics.
Conversely, a congruence constructed in this way has clearly order $1$, because cutting it
with a
$3$--dimensional space we obtain the family of secant lines to a irreducible nondegenerate
curve of degree $3$. Summarizing, we have: 
\begin{teo}\label{unocon}
Every irreducible subvariety of $\Pm^{r+2}$ of degree $3$ and codimension $2$ contains a
$2$--dimensional family of  $(r-1)$--dimensional quadrics, and the family of $\Pm^r$'s
containing those quadrics is a congruence of order $1$ and class $3$, parametrized by a
Veronese surface $\Si$. 
$\Si$ lie on the Schubert cycle $\op_V=\{\s\in G(r,r+2)\,/\,\Pm^r(\s)\supset V\}$, being $V$
the singular locus of all the quadrics, and $r-2\geq\di V\geq r-4$. Conversely, every
congruence of order $1$ with nonlinear irreducible focal locus is constructed in this way.
\qed
\end{teo}
\begin{rem}\label{Veronese}
{\em Thus there are only three cases of congruences in the case II without fixed points, and the
rest are cones over one of those; that is, using the isomorphism 
$$\op_C\cong G(r-\di C-1,r-\di C+1)$$
there are only three different embeddings of the Veronese surface in a
grassmannian as a congruence of bidegree $(1,3)$: in
$G(1,3)$, in
$G(2,4)$ and in
$G(3,5)$.}
\end{rem} 
\section{Case II: $F(\Si)$ is reducible}\label{dos}
In this case, for each $\s\in\Si$, the focal quadric in $\Pm^r(\s)$ is reducible,
$Q(\s)=\Pm_1^{r-1}(\s)\cup\Pm_2^{r-1}(\s)$, and the components of $F(\Si)$ are the
projective realizations, $X$ and $\pi$, of the $\Pm_1^{r-1}(\s)$'s and of the
$\Pm_2^{r-1}(\s)$'s, respectively. By 
(\ref{compo}), $\di X=\di\pi=r$ and $\pi$ is linear.
\begin{propn}\label{scroll}
Under the above assumptions, $X$ is an $r$--dimensional scroll (resp. a curve, for $r=1$)
such that each one of its generators (resp. points, for $r=1$) is contained in exactly one
linear pencil of $\Pm^r$'s of the congruence. 
\end{propn}
{\bf Proof:} For $r=1$, we have already seen that $\Si=\Si(X,\pi)$ in the proof of 
(\ref{compo}). Being
$P\in X$ generic, the lines of the congruence containing $P$ lie on the plane $\pi+P$, hence they
are exactly one linear pencil. 

If $r>1$, the family of $\Pm_1^{r-1}(\s)$'s contained in $X$ cannot be
$2$--dimensional; otherwise $X$ would be a projective space and, as $\Si$ is nondegenerate,
$\di X\cap\pi=(r-2)$, so $X\cap\pi=\Pm_1^{r-1}(\T)\cap\Pm_2^{r-1}(\T)$ for all
$\T\in\Si$; but the family of $\Pm^{r-1}$'s in $X$ containing $X\cap\pi$ is not
$2$--dimensional, a contradiction. Therefore $X$ is a scroll and each
one of its generators $\Pm_1^{r-1}(\s)$ is contained in an infinity of elements of $\Si$, all
of them lying in $\Pm_1^{r-1}(\s)+\pi$; by (\ref{irred}) and (\ref{ratrul}) such
elements are exactly the linear pencil
$\op_{\Pm_1^{r-1}(\s),\Pm_1^{r-1}(\s)+\pi}=\A^{-1}(\A(\s))$. \qed
\begin{example}\label{=1} {\it The case} $r=1$: {\em let $\pi\subset\Pm^3$ be a line,
$X\subset\Pm^3$ an irreducible curve of degree $n$ and suppose $\Si=\Si(\pi,X)$ has order $1$.
By (\ref{irred}), the generic plane $H\supset\pi$ contains exactly one point of $X$
out of $\pi$ (the base point of the pencil of lines of $\Si$ contained in $H$). Thus $\pi$ cuts
$X$ in $n-1$ points counted with multiplicity, and the projection
$$\G:X\lrw\{H\in{\Pm^3}^*\,/\,H\supset\pi\}\cong\Pm^1$$ given by $\G(x)=x+\pi$
shows that $X$ is rational and smooth out of $\pi$. Moreover, if
$H\not\supset\pi$ is generic, it cuts $X$ in $n$ distinct points, no two
collinear with $X\cap\pi$, so $H$ contains exactly $n$ lines of $\Si$; therefore
$\cl(\Si)=\deg(X)=n$.


$X$ is a projection of a rational normal curve
$\Gamma_n\subset\Pm^n$ from a space $V\subset\la x_1,\ldots,x_{n-1}\ra$ 
, $x_1,\ldots,x_{n-1}\in\Gamma_n$, disjoint of $\Gamma_n$. There are two possibilities: either
$X\subset\Pm^3$ is nondegenerate, equivalently $\di V=n-4$, or $X$ is a plane
curve, equivalently $\di V=n-3$. Conversely, if $X$ is such a projection of a rational
normal curve, and $\pi$ is a line containing the projection of $\la
x_1,\ldots,x_{n-1}\ra$, the congruence $\Si(\pi,X)$ has clearly order $1$.\qed}
\end{example}
\begin{rem}\label{=1dos}
{\em The curve $X$ is smooth out of $\pi$, but it is not smooth in general. For instance,
consider the congruence $\Si(X,\pi)$ where $X$ is a nodal cubic in a plane $L$ and $\pi$ is a
line meeting $L$ in the double point of $X$.}
\end{rem}

\begin{example}\label{>1} The case $r>1$: {\em let $\Si\subset G(r,r+2)$ be a congruence of
order $1$ in the case II. In (\ref{scroll}) we have shown how $\Si$ is constructed
from its focal locus $F(\Si)=X\cup\pi$. According to
(\ref{greater}), assume that $\Si$ has no fixed points. By each generator
$\Pm_1^{r-1}(\s)\subset X$ passes the linear pencil
$\op_{\Pm_1^{r-1}(\s),\Pm_1^{r-1}(\s)+\pi}$. Let us denote $\Si(\pi,X)=\Si$.
By (\ref{irred}), the projection
$\G:X\lrw\{H\in{\Pm^{r+2}}^*\,/\,H\supset\pi\}\cong\Pm^1$ given by
$\G(x)=x+\pi$ forces $X$ to be a rational scroll and smooth out of $\pi$. In
order to compute $\cl(\Si)$, take a generic $3$--dimensional subspace
$L\subset\Pm^{r+2}$ and consider the section $\Si_L$, whose
focal locus is $F(\Si_L)=(X\cap L)\cup(\pi\cap L)$. By 
 (\ref{cortar}) and the above example, we get:
$\cl(\Si)=\cl(\Si_L)=\deg(X\cap L)=\deg(X)$. 


Being $\deg X=n$, $X$ is a projection of a rational
normal scroll $R\subset\Pm^{n+(r-1)}$ of the same degree. The center $V$ of such projection
($p_V$) is disjoint of $R$, and is contained in the space generated by $p_V^*(\pi\cap
X)$. Moreover $\di V=(n+(r-1)-\di\la X\ra-1)$, being $\di\la X\ra=r+2$ or $r+1$. 
Conversely, every $r$--dimensional rational normal scroll can be projected to $\Pm^{r+2}$ in
this way, obtaining the focal scroll of a congruence of order $1$ in the case II without
fixed points.\qed} 
\end{example}
The next proposition provides a criterion for the smoothness of a congruence in the case II.
\begin{teo}\label{singular}
Let $\Si=\Si(\pi,X)\subset G(r,r+2)$ be a congruence in the case II. $\Si$ is a rational
ruled surface of degree
$\deg(\Si)=\deg(X)+1$ contained in $\op_{\pi}=\{\s\in G(r,r+2)\,/\,\di(\Pm^r(\s)\cap\pi)\geq
r-1\}$. If
$\Si$ has no fixed points, by $[\pi]$ are passing at most $n-r$ generators of
$\Si$ (exactly $n-1$ if $r=1$) counted with multiplicity. The singular locus
of $\Si$ consists of, at most, $[\pi]$ and the multiple generators of $\Si$ passing by
it. Hence $\Si$ is smooth if and only if $\pi\cap X$ contains at most one generator of
$X$. In particular, if $n\leq r+1$, $\Si$ is always smooth, and if $r=1$, then $\Si$ is smooth iff
$n\leq 2$. 
\end{teo}
{\bf Proof:} If $r=1$, consider
$\G:X\lrw\{H\in{\Pm^3}^*\,/\,H\supset\pi\}\cong\Pm^1$ defined in (\ref{=1}). It
is an isomorphism in $X\setminus\{[\pi]\}$ and puts in correspondence the $n-1$ points
(counted with multiplicity) of $\pi\cap X$ with the planes by $\pi$ containing pencils of
lines of $\Si$ with base point in $\pi\cap X$; equivalently $[\pi]$ belongs to those
pencils. Thus $[\pi]$ will be contained in $n-1$ generators of $\Si$ counted with multiplicity.

If $r>1$ and $\Si$ has no fixed points, consider the projection
$p_V:R\lrw X$ given in (\ref{>1}), whose center $V$ is
contained in $\la p_V^*(\pi\cap X)\ra$. Clearly $p_V^*(\pi\cap X)\sim H-F$
where $\sim$ denotes linear equivalence, $H$ the hyperplane section of
$R$ and $F$ one of its generators. The generators of $R$ contained in $p_V^*(\pi\cap X)$
are in correspondence with the generators of $\Si$ passing by $[\pi]$. Thus, it is enough to
show that $p_V^*(\pi\cap X)$ contains at most $n-r$ generators counted with multiplicity:
a hyperplane section $H\subset R$ is a divisor of the form $H=C+F_1+\ldots+F_k$, where $F_i$ are
generators and
$C=\overline{\bigcup_{F\neq F_1,\ldots,F_k}(F\cap H)}$ is an irreducible
$(r-1)$--dimensional scroll with disjoint generators (since $R$ is
normal); then $\di\la C\ra\geq 2(r-2)+1$ and $\deg(C)\geq\codi(C\subset\la
C\ra)+1\geq r-1$. Since $\deg(H)=n$, $k\leq n-(r-1)=n-r+1$, which is our claim. \qed

Being $n_1,\ldots,n_r\geq 1$ integers, $n=n_1+\ldots+n_r$, $r>1$, let
$R(n_1,\ldots,n_r)\subset\Pm^{n+(r-1)}$ denote the rational normal scroll generated by the
rational normal curves $\Gamma_{n_1},\ldots,\Gamma_{n_r}$ in disjoint spaces (given $r$
isomorphisms $\nu_i:\Pm^1\lrw\Gamma_{n_i}$, $i=1,\ldots,r$, $R(n_1,\ldots,n_r)$ is the scroll
generated by the spaces $\la\nu_1(t),\ldots,\nu_r(t)\ra$, with $t\in\Pm^1$);
$H$ denotes its hyperplane section and $F$ one of its generators.
\begin{propn}\label{lisingu}
The rational normal scroll $R(n_1,\ldots,n_r)$ can be projected to $\Pm^{r+2}$ providing the
focal scroll $X$ of a nonsingular congruence $\Si(\pi,X)$. Hence, given $n\geq r$, there
exist smooth congruences $\Si\subset G(r,r+2)$ of order $1$ and class $n$ without fixed
points. Moreover, if $n\geq r+2$ and if $\max(n_1,\ldots,n_r)\geq 3$, $R(n_1,\ldots,n_r)$
can be projected to $\Pm^{r+2}$ providing the focal scroll of a singular congruence
$\Si(\pi,X)$. Hence, if $n\geq r+2$, there exist singular congruences $\Si\subset G(r,r+2)$
of order $1$ and order $n$.
\end{propn}
{\bf Proof:} As we have seen in the proof of (\ref{singular}), the singularity of
$\Si(\pi,X)$ depends of the number of generators of $R:=R(n_1,\ldots,n_r)$ contained in
$p_V^*(\pi\cap X)\sim H-F$; we will thus have to choose properly a divisor $C\sim H-F$
in each case.

For the first part, it is sufficient to show that $|H-F|$ contains irreducible divisors: the
trace of $|H|$ over $F$ is the complete linear series $|\Te_F(H\cap F)|$, so we have an exact
sequence 
$$
0\rw H^0(\Te_R(H-F))\lrw H^0(\Te_R(H))\lrw H^0(\Te_F(1))\rw 0
$$ 
and $h^0(\Te_R(H-F))=n$; moreover the trace of $|H-F|$ over another generator $F'$ is
$|\Te_{F'}(H\cap F')|$, so the exact sequence
$$
0\rw H^0(\Te_R(H-F-F'))\lrw H^0(\Te_R(H-F))\lrw H^0(\Te_{F'}(1))\rw 0
$$
provides $h^0(\Te_R(H-F-F'))=n-r$. Making $F'$ vary in $R$,
the set of reducible elements of $|H-F|$ has dimension $n-r$, hence it is a proper closed
subset.

For the second part, we will show that $k=\max(n_1,\ldots,n_r)$ is the greatest number of
generators of $R$ contained in a hyperplane (thus, if $k\geq 3$, a divisor in $|H-F|$ can
contain $k-1\geq 2$ generators of $R$): let $F_1,\ldots,F_l$ be generators of $R$, $F_i=\la
P_i^1,\ldots,P_i^r\ra$, $P_i^j\in \Gamma_{n_j}$, $i=1,\ldots,l$. They are contained in a
hyperplane iff $F_1+\ldots+F_l\neq\Pm^{n+r-1}$, iff $\la P_1^1,\ldots,P_l^1\ra+\ldots+\la
P_1^r,\ldots,P_l^r\ra\neq\Pm^{n+r-1}$, iff $l\leq n_j$ for some $j=1,\ldots,r$, iff
$l\leq\max(n_1,\ldots,n_r)$, which is our claim. \qed

The next proposition deals with the existence of congruences with assigned linearly normal
model. We will use the notation of \cite{hartshor} chapter V \S 2.
\begin{propn}\label{exist}
Given $n\geq 1$ and $1\leq r\leq n$, there exist congruences without fixed
points $\Si\subset G(r,r+2)$ in the case II of order $1$ and class $n$, with
given invariant $e$, $0\leq e\leq n-1$, $n-e\equiv 1(\mod 2)$.
\end{propn}
{\bf Proof:} Let $X_e$ be the geometrically ruled surface of invariant $e$,
$C_0$ its minimal directrix ($C_0^2=-e$) and $F$ one of its generators. Suppose $0\leq e\leq
n-1$ and $n-e\equiv 1(\mod 2)$ and consider
the morphisms $$
\begin{array}{cc}
{\psi_1:X_e\lrw\Pm^1,\quad}&{\psi_2:X_e\lrw\Pm^n}\\
\end{array}
$$ 
 given by the linear systems $|F|$
and
$|C_0+(n+e-1)/2\,F|$, respectively. Since $(n-e-1)/2\geq e$ the second one is base
point free, so both maps are regular. Regard the composition 
$$
\psi:X_e\stackrel{\psi_1\times\psi_2}{\lrw}\Pm^1\times\Pm^n\stackrel{i}{\lrw}\op_{\pi}\subset\Pm^{2n+2}
$$
where $\pi\subset\Pm^{n+2}$ is an $n$--dimensional subspace, and we use the
identification of $\op_{\pi}$ with the cone with vertex $[\pi]$ over the Segre
variety given in (\ref{cone}). If $H\subset\Pm^{2n+2}$ is a hyperplane, then 
$\psi^{-1}(H)=C_0+(n+e+1)/2\,F$ and $\psi^*(|H|)=|C_0+(n+e+1)/2\,F|$. Hence
$\psi(X_e)$ is a rational normal ruled surface of degree $n+1$ and invariant
$e$. Since $\psi_1^{-1}(x)\sim F$, the condition of (\ref{ratrul}) holds, so
$\Si$ has order $1$. Each
generator $F\subset\psi(X_e)$ parametrizes a pencil of $\Pm^n$'s not containing
$\pi$, equivalently, their focal locus does not lie on $\pi$, so $\psi(X_e)$ is
in the case II, and it is smooth since $[\pi]\notin\psi(X_e)$. 

Moreover, for $r<n$ it is sufficient to take the
generic projections of $\psi(X_e)$ according to 
(\ref{cortar}).\qed 
\begin{rem}\label{expect}{\em Moreover, every congruence in the case II is
expected to be a projection of one of those constructed above,
because, as rational surfaces, those are their linearly normal models.}
\end{rem}
\section{Case III: $F(\Si)$ is an $r$-dimensional projective space}\label{tres}
$F(\Si)$ is the only $r$-dimensional projective
subspace
$\pi\subset\Pm^{r+2}$ such that 
$\Si\subset\op_{\pi}=\{\s\in G(r,r+2)\,/\,\di\Pm^r(\s)\cap\pi\geq r-1\}$. Moreover, in this
case,
$2(\Pm^{r-1}(\s)):=2(\Pm^r(\s)\cap\pi)$ is the focal quadric of $\Si$ at
$\Pm^r(\s)$. The fibers of the map $$\A:\Si\setminus\{[\pi]\}\lrw
H_{\pi}\cong\Pm^1$$ defined in (\ref{ratrul}) parametrize pencils whose 
base loci lie in $F(\Si)=\pi$, so they all contain $[\pi]$. Therefore {\it $\Si$ is
a rational cone with vertex $[\pi]$}.
Consider the projection $\B:\Si\setminus\{[\pi]\}\lrw\pi^*$ defined in 
(\ref{fixed}), $\,\B(\s)=\Pm^r(\s)\cap\pi$: since the fibers of this map are
$\B^{-1}(\Pm^{r-1}(\s))=\A^{-1}(\Pm^r(\s)+H)$, $\B(\Si)$ is a curve. Thus we have
a birrational map 
\begin{equation}
\setlength{\unitlength}{5mm}
\begin{picture}(7,2.5)
\put(0,2){\makebox(0,0){$H_{\pi}$}}
\put(0,0){\makebox(0,0){$H$}}
\put(7,2){\makebox(0,0){$\B(\Si)\subset\pi^*$}}
\put(7,0){\makebox(0,0){$\bigcap_{\s\in\A^{-1}(H)}\Pm^r(\s)$}}
\put(1,2){\vector(1,0){4}}
\put(1,0){\vector(1,0){2.5}}
\put(1,0.1){\line(0,-1){0.2}}
\put(3,2.5){\makebox(0,0){$\Lambda$}}
\end{picture}
\label{Lambdada}
\end{equation}  
such that $\Lambda\circ\A=\B$ and so $\Si$ is characterized by $\Lambda$ in the
next way: 
\begin{equation}
\Si=\bigcup_{H\supset\pi}\op_{\Lambda(H),H}
\label{caracteriza11}
\end{equation} 
being $\op_{\Lambda(H),H}=\{\s\in G(r,r+2)\,/\,\Lambda(H)\subset\Pm^r(\s)\subset H\}$.

$\B(\Si)$ parametrize a family of hyperplanes in $\pi$, hence it is
developable, and we can say how it is constructed: suppose $\Si$ has no fixed points (otherwise
see (\ref{greater})); 
$\B(\s)$ is the family of
$(r-1)$--osculating hyperplanes to a nondegenerate curve $C\in\pi$, which is
birrational to
$\B(\Si)$. Thus we have another birrational map:
\begin{center}
\setlength{\unitlength}{5mm}
\begin{picture}(10.2,2.5)
\put(-3,2){\makebox(0,0){$H_{\pi}$}}
\put(-3,0){\makebox(0,0){$H$}}
\put(3,2){\makebox(0,0){$C\subset\pi$}}
\put(9.3,0){\makebox(0,0){$\Lambda'(H)=(r-1)$--th focal locus of
$\B(\Si)$ at 
$\Lambda(H)$}}
\put(-2,2){\vector(1,0){3}}
\put(-2,0){\vector(1,0){2.5}}
\put(-2,0.1){\line(0,-1){0.2}}
\put(-0.5,2.5){\makebox(0,0){$\Lambda'$}}
\end{picture}
\end{center}  
that characterizes:
\begin{equation}
\Si=\bigcup_{H\supset\pi}\op_{T_{r-1,C,\Lambda'(H)},H}
\label{caracteriza21}
\end{equation}
where $T_{r-1,C,\Lambda'(H)}$ denotes the $(r-1)$-osculating hyperplane to $C$ at
$\Lambda'(H)$.
 
Let $\phi$ ($=\A\times\B$) be the map defined in (\ref{cone}); $\Si$ is a cone with vertex
$[\pi]$ over the curve
$\phi(\Si\setminus\{[\pi]\})\subset\pi^*
\times H_{\pi}$, which is a hyperplane section of $\Si$. Since this curve is the graph of the
map $\Lambda:H_{\pi}\lrw\pi^*$, $\phi(\Si\setminus\{[\pi]\})$ is an smooth curve of
degree $n+1$ (where 
$n:=\deg(\B(\Si))=
\deg(\Lambda(H_{\pi}))$). Thus $\deg(\Si)=\deg(\phi(\Si\setminus\{[\pi]\}))
=n+1=\deg(\Lambda(H_{\pi}))+1$ and we conclude $\cl(\Si)=n$. Moreover the singular locus of
$\Si$ is only the vertex 
$[\pi]$. Summarizing, we have:
\begin{teo}\label{trescon1}
For every congruence $\Si\subset G(r,r+2)$ in the case III without fixed points there exist a
regular map
$$\Lambda:H_{\pi}\lrw\pi^*$$ (being $\pi=F(\Si)$ and $H_\pi$ the pencil of hyperplanes
containing $\pi$) with nondegenerate image
$\Lambda(H_{\pi})$ such that: 
$$
\Si=\bigcup_{H\supset\pi}\op_{\Lambda(H),H}
$$
being $\op_{\Lambda(H),H}=\{\s\in
G(r,r+2)\,/\,\Lambda(H)\subset\Pm^r(\s)\subset H\}$. Moreover, $\Si$ is a cone of degree 
$\deg(\Lambda(H_{\pi}))+1$ ($\cl(\Si)=n$), with vertex $[\pi]$ over an smooth curve.
Conversely, a congruence constructed in this way is in the case III. \qed
\end{teo}


The curve $\B(\Si)\subset\pi^*$ is the projection of its normal model
$\Gamma_n\subset(\Pm^n)^*$  from a subspace $V\subset(\Pm^n)^*$, hence the map
$\Lambda$ factorizes through $\Gamma_n$. Regard   
$\pi\subset\Pm^n$, and everything contained in a projective space $\Pm^{n+2}$
such that $\pi$ is obtained cutting $\Pm^n$ with the 
$(r+2)$--dimensional space containing the congruence $\Si$ (allowing us to identify the
hyperplanes in $\Pm^{r+2}$ containing $\pi$ with the hyperplanes in $\Pm^{n+2}$ containing
$\Pm^{n}$), we will have a map 
$\Lambda_n:H_{P^n}\lrw(\Pm^n)^*$ providing a congruence $\Si_n\subset G(n,\Pm^{n+2})$ whose
projection to $G(r,\Pm^{r+2})$ is $\Si$. Hence we have a commutative diagram 
\begin{center}
\setlength{\unitlength}{5mm}
\begin{picture}(18.2,3.5)
\put(1.5,3){\makebox(0,0){$\Gamma_n\subset(\Pm^n)^*$}}
\put(1,0){\makebox(0,0){$\B(\Si)\subset\pi^*$}}
\put(9,3){\makebox(0,0){$(\Pm^n)^*\times H_{P^n}$}}
\put(9,0){\makebox(0,0){$(\pi)^*\times H_{\pi}$}}
\put(16,3){\makebox(0,0){$\op_{P^n}\supset\Si_n$}}
\put(16,0){\makebox(0,0){$\op_{\pi}\supset\Si$}}
\put(6,3){\vector(-1,0){2}}
\put(14,3){\vector(-1,0){2}}
\put(6,0){\vector(-1,0){2}}
\put(14,0){\vector(-1,0){2}}
\put(0,2.2){\vector(0,-1){1.5}}
\put(2.3,2.2){\vector(0,-1){1.5}}
\put(9,2.2){\vector(0,-1){1.5}}
\put(15,2.2){\vector(0,-1){1.5}}
\put(17,2.2){\vector(0,-1){1.5}}
\put(13,3.5){\makebox(0,0){$\phi_n$}}
\put(13,0.5){\makebox(0,0){$\phi$}}
\put(18.2,1.45){\makebox(0,0){$\rho_{\Pm^{r+2}}$}}
\end{picture}
\end{center}  
We have thus proved:
\begin{teo}\label{trescon}
Being $n\geq 1$, there exist a congruence $\Si_n\subset G(n,n+2)$ constructed in
the next way: given $L\subset\Pm^{n+2}$ an $n$--dimensional projective space
and $\Lambda:\Pm^1\cong H_L\lrw L^*$ a $n$--th Veronese embedding,
$\Si_n:=\bigcup_{H\supset L}\op_{\Lambda(H),H}$. $\Si_n$ is a  rational normal cone of
degree $n+1$. Every congruence $\Si\subset G(r,r+2)$ of order $1$ and class $n$ without fixed
points in the case III can be obtained cutting $\Si_n$ with a suitable projective
subspace $\Pm^{r+2}$, that is:
$$
\Si_n\!\subset\! G(n,n+2)\stackrel{\rho_{\Pm^{r+2}}}{\lrw}\Si\!\subset\!
G(r,r+2).
$$\qed
\end{teo}
\section{Smoothness of a congruence of order $1$}\label{lisas}
According to (\ref{greater}), in order to study the smoothness of congruences of order
$1$ we only need to consider congruences without fixed points, translating the results to the
general case. Summarizing 
(\ref{unocon}),
(\ref{singular}), (\ref{lisingu}) and
(\ref{trescon1}), we have the next theorems: 
\begin{teo}\label{conclusiones}
Let $r,\,s$ be two integers such that $r\geq 1$ and $-1\leq s\leq r-1$. If $\Si\subset G(r,r+2)$
is a surface of order $1$ and class $n$ such that $\di(T(\Si))=s$, then 
$n\geq r-s-1$. Furthermore:
\begin{itemize}
\item If $s=r-1$, $\Si$ is a plane.
\item For every $s\leq r-2$, there exist singular surfaces $\Si\subset G(r,r+2)$ of bidegree 
$(1,n)$ such that $\di(T(\Si))=s$.
\end{itemize}
If $s=r-2$, then:
\begin{enumerate}
\item If $n=1,2$ and $\Si$ is not a cone, $\Si$ is smooth.
\item If $n\geq 4$, $\Si$ is singular.
\item If $n=3$, there exist smooth and singular surfaces $\Si\subset G(r,r+2)$ of
bidegree $(1,3)$ such that $\di T(\Si)=s$.
\end{enumerate}
If $-1\leq s\leq r-3$, then:
\begin{itemize}
\item If $n=r-s-1$ or $n=r-s$ and $\Si$ is not a cone, $\Si$ is smooth.
\item If $n\geq r-s+1$, there exist smooth and singular surfaces 
$\Si\subset G(r,r+2)$ of bidegree $(1,n)$ such that $\di T(\Si)=s$. 
\end{itemize} 
\end{teo}
{\bf Proof:} Let $\Si\subset G(r,r+2)$ be a congruence of order $1$ and class $n$ such that 
the dimension of its fixed locus is $\di(T(\Si))=s$; if $s=r-1$, $\Si$ is the family of 
$\Pm^r$'s containing $T(\Si)$, which is parametrized by a plane; in this case  
$n=0=r-s-1$; suppose now that $s\leq r-2$. Applying 
(\ref{greater}), $\Si$ can be projected isomorphically to a congruence $\Si'\subset
G(r-s-1,r-s+1)$ without fixed points, having order $1$ and class $n$; if $\Si'$ is in the case 
I, then $n=3$ and $r-s-1\in\{1,2,3\}$, and then $n\geq r-s-1$; otherwise 
(\ref{fixed}) shows that 
$n=\cl(\Si')\geq r-s-1$. Moreover for every $n\geq r-s-1$, there exist a congruence 
$\Si'\subset G(r-s-1,r-s+1)$ of class 
$n$ in the case III without fixed points (see (\ref{trescon})), which is 
parametrized by a cone. Excepting those congruences, the rest of the theorem is consequence
of (\ref{singular}) and (\ref{unocon}) and (\ref{lisingu}). \qed

\begin{teo}\label{lisascon}
Being $r\geq 1$, the only smooth surfaces in $G(r,r+2)$ of bidegree $(1,n)$ are:
\begin{enumerate}
\item The plane parametrizing the family of $\Pm^r$'s containing a fixed $\Pm^{r-1}$  
($n=0$).
\item The Veronese surface, described in 
(\ref{unocon}), for $r=1,2,3$ ($n=3$); embedding $G(1,3)$, $G(2,4)$ or 
$G(3,5)$ through the map given in (\ref{greater}), it can be consider in every $G(r,r+2)$.
\item The rational ruled surfaces $\Si(\pi,X)$ described in Section \ref{dos} ($n=\deg X$), and 
verifying that $\pi$ contains at most one generator of $X$; this condition is always verified for
$n=r,r+1$.
\item The surfaces described in 2. and 3. embedded in a higher dimension
grassmannian $G(r,r+2)$ through the isomorphism exposed in (\ref{greater}). \qed
\end{enumerate}
\end{teo}
Let us finally enumerate the smooth congruences of order $1$ in $G(1,3)$, recovering a Ziv Ran's
result given in \cite{ran}. We remark that, for a congruence $\Si(\pi,X)$ in
$G(1,3)$, $\pi$ is a line cutting the curve $X$ in $n-1$ points; hence, imposing
the smoothness condition of (\ref{singular}), we have $n-1\leq 1$, that is $\deg X\leq
2$. 
\begin{cor}\label{zivran}
The only smooth congruences in $G(1,3)$ of order $1$ are:
\begin{enumerate}
\item The plane parametrizing the lines passing by a fixed point (bidegree (1,0)).
\item The Veronese surface $\Si(C)$ parametrizing the secant lines to a rational normal cubic
$C$ (bidegree (1,3)).
\item The quadric $\Si(\pi,\pi')$ parametrizing the lines cutting two given lines $\pi$ and $\pi'$
(bidegree (1,1)).
\item The rational normal cubic $\Si(\pi,C)$ parametrizing the lines cutting a given line $\pi$
and a given irreducible conic $C$ meeting in a point $P=\pi\cap C$ (bidegree (1,2)). \qed
\end{enumerate}
\end{cor}

\begin{quote}
\emph{Authors' address:} Departamento de Algebra,
Universidad de Santiago de Compostela, 15706 Santiago de Compostela,
Galicia, Spain. Phone: 34-81563100-ext.13152. Fax: 34-81597054. {\tt e}-mail: {\tt
pedreira@zmat.usc.es}
\end{quote} 
\end{document}